\documentclass[11pt]{amsart}

\newcommand{\eq}{~\eqref}
\numberwithin{equation}{section}

\newtheorem{thm}{Theorem}[section]
\newtheorem{mthm}{Theorem}
\newtheorem{lem}[thm]{Lemma}

\newtheorem{pro}[thm]{Proposition}

\theoremstyle{definition}

\theoremstyle{remark}

\newtheorem{rem}[thm]{Remark}

\newcommand{\ra}{{\mbox{$\rightarrow$}}}
\newcommand{\R}{{\mathbf{R}}}                     %% Real Numbers
\newcommand{\s}{{\mathbf{S}}}                     %% Sphere
\newcommand{\W}{\mathcal{D}}                %% Differences of support fcns
\newcommand{\cn}{\colon}               %% Colon
\catcode`@=11 \@mparswitchfalse  %This puts the \mnote's on the
\newcounter{mnotecount}[section]

\title{A Generalized Affine Isoperimetric Inequality}

\author[Chen]{Wenxiong Chen${}^*$}
\thanks{${}^*$ Chen supported, in part, by NSF Grant DMS-0072328}
\address{Department of Mathematics,
Yeshiva University, New York NY 10033, USA} \email{wchen@yu.edu}
\urladdr{http://math.smsu.edu/$\sim$wchen/}

\author[Howard]{Ralph Howard${}^\dagger$}
\thanks{${}^\dagger$ Howard partially supported by ONR-DEPSCoR
Contract \# N000140310675} 
\address{Department of Mathematics, University of South
Carolina,
Columbia, SC 29208, USA}
\email{howard@math.sc.edu}
\urladdr{www.math.sc.edu/$\sim$howard}

\author[Lutwak]{Erwin Lutwak${}^\ddagger$}
\thanks{${}^\ddagger$ Lutwak, Yang, and Zhang supported, in part, by NSF Grant DMS-0104363}
\address{Department of Mathematics,
Polytechnic University,
Brooklyn, NY 11201, USA}
\email{elutwak@duke.poly.edu}
\author[Yang]{Deane Yang${}^\ddagger$}
\address{Department of Mathematics,
Polytechnic University,
Brooklyn, NY 11201, USA}
\email{dyang@duke.poly.edu}
\urladdr{http://www.math.poly.edu/$\sim$yang/}

\author[Zhang]{Gaoyong Zhang${}^\ddagger$}
\address{Department of Mathematics,
Polytechnic University,
Brooklyn, NY 11201, USA}
\email{gzhang@duke.poly.edu}
\urladdr{http://www.math.poly.edu/$\sim$gzhang/}

\date{\today}

\begin{document}
\allowdisplaybreaks

\begin{abstract}
A purely analytic proof is given for an inequality that has as a
direct consequence the two most important affine isoperimetric
inequalities of plane convex geometry: The Blaschke-Santalo inequality
and the affine isoperimetric inequality of affine differential
geometry.
\end{abstract}

\maketitle

%%%%%%%%%%%%%%%%%%%%%%%%%%%%%%%%%%%%%%%%%%%%%%%%%%%%%%%%%%%%%%%%%%%%%%
\section{Introduction.}\label{sec:intro}
%%%%%%%%%%%%%%%%%%%%%%%%%%%%%%%%%%%%%%%%%%%%%%%%%%%%%%%%%%%%%%%%%%%%%%

In \cite{Harrell:direct}, Harrell showed how an analytic approach
could be used to obtain a well-known \emph{Euclidean} inequality
of plane convex geometry -- the Blaschke-Lebesgue inequality. In
this article we show how a purely analytic approach can be used to
establish the best known \emph{affine} inequalities of plane
convex geometry. To be precise, we will use a purely analytic
approach to establish an analytic inequality that has as an
immediately consequence both the affine isoperimetric inequality
of affine differential geometry and the Blaschke-Santal\'o
inequality. What's more significant, we are able to remove the
\lq\lq convexity\rq\rq\ assumption and thus establish an
inequality with applications to the planar $L_p$ Minkowski problem
(see, e.g., \cite{Lutwak:B-M-F}, \cite{Lutwak-Oliker}, \cite{St1},
\cite{St2}, \cite{U}) with not necessarily positive data.
%\bigpagebreak

Let $C\subset \R^2$ be a compact convex set. Let $\s$ be the unit
circle parameterized by
%the angular coordinate $\theta$ (corresponding
%to $(\cos\theta,\sin\theta)\in\s$).  We will use the notation
$$
e(\theta):=(\cos\theta,\sin\theta).
$$
Then $h=h_C\cn \s\to
\R$ defined by
$$
h(\theta) :=\max_{x\in C} e(\theta)\cdot x
$$
is the \emph{support function} of $C$.

The affine isoperimetric inequality of affine differential
geometry states that if a plane convex figure has support function
$h\in C^2(\s)$, then
\begin{equation}\label{1}
\gathered
  4\pi^2
  \int_\s
  h(h+h'')\,d\theta
  \ge
  \bigg(
  \int_\s
  (h+h'')^{2/3}\,d\theta
  \bigg)^3\\
\text{with equality if and only if the figure is an ellipse.}
\endgathered
\end{equation}
The integral on the left is twice the area of the figure, while the
integral on the right is the so called \emph{affine perimeter} of the
figure.

The Blaschke-Santal\'o inequality
states that if a convex figure is \emph{positioned} so that its support
function $h$ is positive and
\begin{equation}\label{2a}
  \int_\s \frac{\sin\theta}{h(\theta)^3} \, d\theta
  =
  0
  =
  \int_\s\frac{\cos\theta}{ h(\theta)^{3}} \, d\theta,
\end{equation}
then
\begin{equation}\label{2}
\gathered
  4\pi^2\bigg(
  \int_\s h^{-2} d\theta
  \bigg)^{-1}
  \ge
  \int_\s h(h+h'')\, d\theta,\\
\text{with equality if and only if the figure is an ellipse.}
\endgathered
\end{equation}
The integral on the left is twice the area of the \emph{polar reciprocal}
of the figure.
%The H\"older inequality tells us that
%the Blaschke-Santal\'o inequality
%is stronger (i.e., directly implies) the classical isoperimetric
%inequality between perimeter and area. To see this note that the
%H\"older inequality shows that
%$$\gathered
%  \bigg(
%  \int_\s h\, d\theta
 % \bigg)^{2}
  %\ge
  %8\pi^3
  %\bigg(
  %\int_\s h^{-2} d\theta
  %\bigg)^{-1}\\
%\text{with equality if and only if the figure is a circle.}
%\endgathered
%$$ The integral on the left is the perimeter of the figure. When
%this inequality is combined with \eq{2} the result is the
%classical isoperimetric inequality in the plane.

In \cite{Lutwak:Bl-Sa}, it was shown that both inequalities \eq{1} and
\eq{2} are encoded in the following inequality: If $K$ and $L$ are
convex figures whose support functions are such that $h_L\in C^2(\s)$
and $h_K$ arbitrary then
\begin{equation}\label{3}
\gathered
  4\pi^2
  \bigg(
  \int_\s
   (h_L+h_L'')h_K\,d\theta
  \bigg)^2
  \ge
  \bigg(
  \int_\s
  (h_L+h_L'')^{2/3} d\theta
  \bigg)^3
  \int_\s
  [h_K^2-(h_K')^2]\,d\theta\\
\text{with equality if and only if $K$ and $L$ are homothetic ellipsoids.}
\endgathered
\end{equation}

Using this version of the inequality, we see that choosing $L=K$
in \eq{3} immediately gives \eq{1}. To see how \eq{3} gives \eq{2}
choose the figure $L$ so that $h_L$ satisfies the equation $h_L''+
h_L =h_K^{-3}$.
%(\mnote{\ralph A little more detail added.  The old
%version is still in the file, commented out}
%This is always possible as
%$h_K$ satisfies \eq{2a} and therefore the equation $y''+y=h_K^{-3}$
%can be solved for a $2\pi$ periodic function, for example by the
%method of Fourier series.  Then $h_L:=y$ will be the support function
%of a convex set by Proposition~\ref{spt-char} below.)
%
%(Note that since $h_K$ satisfies \eq{2a} this is
%always possible.)

In \cite{Lutwak:Bl-Sa}, it was shown that \eq{3} is a consequence
of \eq{1} and the \emph{mixed area} inequality. The aim of this
paper is to establish an analytic inequality that extends
inequality \eq{3}. Our proof of this new analytic inequality uses
none of the tools of convex geometry. We are thus able to remove
the \lq\lq convexity\rq\rq\ assumption, $h''+h\ge 0$, from all the
inequalities above.

%%%%%%%%%%%%%%%%%%%%%%%%%%%%%%%%%%%%%%%%%%%%%%%%%%%%%%%%%%%%%%%%%%%%%%
\section{The main inequality.}\label{sec:ineq}
%%%%%%%%%%%%%%%%%%%%%%%%%%%%%%%%%%%%%%%%%%%%%%%%%%%%%%%%%%%%%%%%%%%%%%

Let $H^1(\s)$ be the Hilbert space of functions $u\cn \s\to \R$
%such that $u$ is absolutely continuous and $u'\in L^2(\s)$.  We
%use
with the norm
$$
\|u\|_{H^1}=\left(\int_\s [u^2+(u')^2]\,d\theta\right)^{\frac12}.
$$
%The space $H^1(\s)$ can also be described as the space of functions
%whose first distributional derivative is in $L^2$.  The norm is a
%Hilbert space norm with corresponding inner product $\langle
%u,v\rangle_{H^1}=\int_\s (uv+u'v')\,d\theta$.

\begin{mthm}[Two Dimensional Analytic Affine Isoperimetric Inequality]
\label{thm:main}\null\

\noindent
Assume

i) $F$ and $h$ are non-negative $2\pi$ periodic functions that do not
 vanish identically.

ii) $F$ is integrable on $S$ and satisfies the orthogonality
conditions
\begin{equation}\label{1.0}
\int_\s F(\theta ) \cos \theta  \, d\theta  = 0 = \int_\s F(\theta
) \sin \theta \, d\theta  .
\end{equation}

iii) $h\in H^1(\s)$.

Then
\begin{equation}\label{1.2}
\left(\int_\s F(\theta ) h(\theta ) \, d\theta \right)^2 \geq
\frac{1}{4 \pi^2} \left(\int_\s F^{2/3} \, d\theta \right)^3
\left(\int_\s [h^2-(h')^2] \, d\theta \right).
\end{equation}
Equality holds if and only if there exist $k_1, k_2, a>0$, and
$\alpha\in \R$ such that
\begin{equation}
h(\theta )=k_1 \sqrt{a^2\cos^2 (\theta -\alpha) + a^{-2}\sin^2 (\theta -\alpha)}
\label{1.3}
\end{equation}
and $F$ is given almost everywhere by
\begin{equation}\label{1.4}
F(\theta )=k_2 ( a^2\cos^2 (\theta -\alpha) + a^{-2}\sin^2 (\theta
-\alpha) )^{-3/2}.
\end{equation}
\label{mthm}
\end{mthm}

\begin{rem}\label{spt-elp}
The functions $h(\theta)$ of the form  \eq{1.3} are exactly
support functions of the ellipses centered at the origin.
\end{rem}

The main ingredient in the proof of Theorem~\ref{thm:main} is a
family of transforms, that leave a few key integrals invariant and
which let us construct maximizing sequences. We will introduce the
transforms and study their properties in Section~\ref{sec:trans}.  In
Section~\ref{sec:proof}, we prove the inequality. In
Sections~\ref{sec:spts-fcns} we study some regularity results for
support functions of planar convex sets.  These regularity results are
used in Section~\ref{sec:aff-iso} to derive the affine isoperimetric
inequality for general planar sets.

%%%%%%%%%%%%%%%%%%%%%%%%%%%%%%%%%%%%%%%%%%%%%%%%%%%%%%%%%%%%%%%%%%%%%%
\section{Function spaces associated with the inequality.}\label{sec:spts-fcns}
%%%%%%%%%%%%%%%%%%%%%%%%%%%%%%%%%%%%%%%%%%%%%%%%%%%%%%%%%%%%%%%%%%%%%%

Because of the integral $\int_\s [h^2-(h')^2]\,d\theta$ that
appears in Theorem~\ref{thm:main}, the natural function space for
the functions $h$ in the theorem is $H^1(\s)$.  Moreover in the
geometric applications Theorem~\ref{thm:main} a natural choice for
the function $h$ is to be a support function of a bounded convex
set and $H^1(\s)$ contains all the support functions.
%However,
%for the geometric applications mentioned in the
%\mnote{\ralph A bit of rewording.}
%introduction, the function $F$ is taken to be
%$F=(h+h'')^2$ where $h$ is a support function.  For a general
%support function the second derivative only exists in a
%generalized sense, say as a distribution, and therefore the
%expression $(h+h'')^2$ is not necessarily defined. In fact for a
%figure as simple as a polygon the $h+h''$ is a sum of point masses
%(i.e. delta functions) and so $(h+h'')^2$ is undefined.
The
following characterizes the support functions of bounded convex
sets.

\begin{pro}\label{spt-char}
A continuous
%\mnote{\ralph Added ``continuous''}
function $h\cn \s\to \R$
%\mnote{\ralph Changed ``$h\cn \R\to \R$'' to ``$h\cn \s\to \R$''.}
is the support function of a bounded convex set if and only if the
second distributional derivative $h''$ satisfies $(h''+h)\ge0$ as a
distribution.
\end{pro}

The proof of this proposition is elementary, and is left to the
reader.

%As a distribution is positive if and only if it is represented by a
%non-negative measure we see that the distribution $h''+h$ is
%represented by a non-negative measure if and only if $h$ is a support
%function.
We now describe the smallest function space that contains the
support functions of convex sets.  Let $\W$ be the set of $2\pi$
periodic functions $u$ such that the distributional derivative
$u''$ is a signed measure.
%(by Riesz's characterization of the
%dual of $C(\s)$ as space of signed measures this is the same as
%$u''$ being a continuous linear functional on $C(\s)$).
The total
variation $\|\mu\|_{TV}$ of a signed measure $\mu$ on $\s$ is its
norm as a linear functional on $C(\s)$.  That is
$$
\|\mu\|_{TV}:=\sup\left\{\int_\s \phi(\theta)\,d\mu(\theta) : \phi\in
C(\s), |\phi(\theta)|\le 1\right\}.
$$
The standard norm on $\W$ is $\|h\|_{L^\infty}+\|h''\|_{TV}$, but, for
geometric reasons, we use the equivalent
norm
$$
\|h\|_\W:=\|h\|_{L^\infty}+\|h''+h\|_{TV}.
$$
The space $\W$ can also be defined as the functions $h$ on $\s$ that
are absolutely continuous and such that the first derivative $h'$ is
of bounded variation.  As functions of bounded variation are bounded
this implies all elements of $\W$ are Lipschitz.  Therefore the
imbedding $\W\subset C^\alpha(\s)$ is compact for $\alpha\in
[0,1)$,
%\mnote{\ralph Added definition of $C^\alpha((\s)$}
where $C^\alpha(\s)$ is the space of H\"older continuous functions
%$u$ such that the norm
%$\|u\|_{C^\alpha}:=\sup_{\theta\in\s}|u(\theta)|
%+\sup_{\theta_1\ne\theta_2}|\theta_2-\theta_1|^{-\alpha}|u(\theta_2)-u(\theta_1)|$
%is finite.

\begin{thm}\label{con-diff}
The space $\W(\s)$ contains all the support functions of bounded
convex sets. Moreover every element of $\W(\s )$ is a difference of
two support functions. Thus $\W(\s)$ is the smallest function space
containing all the support functions.  More precisely if $f\in\W$
there are support functions $h_1,h_2\in \W$ with $f=h_1-h_2$ and
\begin{equation}\label{f-diff}
\|h_1''+h_1\|_{TV},\ \|h_2''+h_2\|_{TV}\le 3\|f''+f\|_{TV}.
\end{equation}
\end{thm}

\begin{proof}  We have already seen that $\W$ contains all the
support functions of bounded convex sets.  Let $f\in \W(\s )$.  Then
$f''+f$ is a signed measure. We now claim that we can write
$f''+f=\mu_+-\mu_-$ where $\mu_+$ and $\mu_-$ are
non-negative measures with the extra conditions that
\begin{equation}\label{no-per}
\int_\s \cos\theta\,d\mu_+=
\int_\s \cos\theta\,d\mu_-=
\int_\s \sin\theta\,d\mu_+=\int_\s \sin\theta\,d\mu_-=0
\end{equation}
and
\begin{equation}\label{mu-bds}
\|\mu_+\|_{TV},\ \|\mu_-\|_{TV}\le 3\|f''+f\|_{TV}.
\end{equation}
To start let $f''+f=\nu_+-\nu_-$ be the Jordan
decomposition (cf.~\cite[p.~274]{Royden:analysis}) of
$f''+f$.  Then $\nu_+$ and $\nu_-$ are non-negative measures and
$\|f''+f\|_{TV}=\|\nu_+\|_{TV}+\|\nu_-\|_{TV}$.  From the definition
of the second distributional derivative (which is formally just
integration by parts)
$$
\int_\s(f''+f)\cos\theta\,d\theta
=-\int_\s f\cos\theta\,d\theta+\int_\s f\cos\theta\,d\theta=0
$$
and likewise $\int_\s(f''+f)\sin\theta\,d\theta=0$.
Using this in $f''+f=\nu_+-\nu_-$ gives
$$
\int_\s \cos\theta\,d\nu_+=\int_\s \cos\theta\,d\nu_-,\qquad
\int_\s \sin\theta\,d\nu_+=\int_\s \sin\theta\,d\nu_-.
$$
Set
\begin{align*}
a&:=\frac1\pi\int_\s \cos\theta\,d\nu_+=\frac1\pi\int_\s \cos\theta\,d\nu_-,\\
b&:=\frac1\pi\int_\s \sin\theta\,d\nu_+=\frac1\pi\int_\s \sin\theta\,d\nu_-.
\end{align*}
Let $C>0$, to be chosen shortly, and set
\begin{align*}
\mu_+&=\nu_++(C-a\cos\theta-b\sin\theta)\,d\theta\\
\mu_-&=\nu_-+(C-a\cos\theta-b\sin\theta)\,d\theta.
\end{align*}
There is an $\alpha$ so that
$a\cos\theta+b\sin\theta=\sqrt{a^2+b^2}\,\cos(\theta +\alpha)$.
Thus if $C:=\sqrt{a^2+b^2}$ the measures $\mu_+$ and $\mu_-$ are
non-negative. Using that in $L^2$ the function $\cos\theta$ is
orthogonal to $\sin\theta$ and to the constants and that
$\int_\s\cos^2\theta\,d\theta=\pi$
\begin{align*}
\int_\s \cos\theta\,d\mu_+&=
\int_\s
\cos\theta\,d\nu_++\int_\s\cos\theta(C-a\cos\theta-b\sin\theta)\,d\theta\\
&=\int_\s \cos\theta\,d\nu_+-a\int_\s\cos^2\theta\,d\theta=0
\end{align*}
and likewise all the other conditions of\eq{no-per} hold.
As $\|f''+f\|_{TV}=\|\nu_+\|_{TV}+\|\nu_-\|_{TV}$ we have
$\min\{\|\nu_+\|_{TV},\|\nu_-\|_{TV}\}\le \frac12\|f''+f\|_{TV}$.
The formulas for $a$ imply
\begin{align*}
|a|&\le\min\left\{\frac1\pi\int_\s |\cos\theta|\,d\nu_+(\theta),\
\frac1\pi\int_\s |\cos\theta|\,d\nu_-(\theta)\right\}\\
&\le
\frac1\pi \min\big\{\|\nu_+\|_{TV},\ \|\nu_-\|_{TV}\big\}\le
\frac1{2\pi}\|f''+f\|_{TV}.
\end{align*}
Likewise $|b|\le (2\pi)^{-1}\|f''+f\|_{TV}$. Using this in the
definition of $C$ gives $ C\le {\sqrt{2}}(2\pi)^{-1}\|f''+f\|_{TV}$.
As $\mu_+$ and $\mu_-$ are non-negative measures their total variation
is just their total mass.  Thus
\begin{align*}
\|\mu_\pm\|_{TV}&=\int_\s 1\,d\mu_\pm=\int_\s 1\,\nu_\pm
+\int_\s(C-a\cos\theta-b\sin\theta)\,d\theta\\
&=\|\nu_{\pm}\|+2\pi C\le (1+\sqrt{2})\|f''+f\|_{TV}\\
&\le 3 \|f''+f\|_{TV}
\end{align*}
This shows that \eq{mu-bds} holds.

We  claim that there is a function $h_+$ so that
$h_+''+h_+=\mu_+$.  To see this expand $\mu_+$ in a Fourier series and
use the equations\eq{no-per} to see that the coefficients of
$\sin$ and $\cos$ vanish.
$$
\mu_+=\frac{a_0}2+\sum_{k=2}^\infty (a_k\cos(k\theta)+b_k\sin(k\theta)).
$$
Then $h_+$ is given explicitly by
$$
h_+(\theta)=\frac{a_0}2+\sum_{k=2}^\infty
\frac{a_k\cos(k\theta)+b_k\sin(k\theta)}{1-k^2}.
$$
The
%\mnote{\ralph Added argument that $h_+$ is continuous.}
formulas $a_k=\pi^{-1}\int_{\s}\cos(k\theta)\,d\mu_+(\theta)$,
$b_k=\pi^{-1}\int_{\s}\sin(k\theta)\,d\mu_+(\theta)$ imply that
$|a_k|, |b_k|\le \frac12 \|\mu_+\|_{TV}$.  Therefore the series
defining $h_+$ converges uniformly and thus $h_+$ is continuous.
Likewise there is a continuous function $h_-$ with
$h_-''+h_-=\mu_-$. As $\mu_+$ and $\mu_-$ are non-negative
measures and formal differentiation of Fourier series corresponds
to taking distributional derivatives, both $h_+$ and $h_-$ are
support functions.

Let $y=f-(h_+-h_-)$.  Then $y''+y=0$.  This implies
$y=\alpha\cos\theta+\beta\sin\theta$ for some constants $\alpha$
and $\beta$.  Thus $f=(h_++\alpha\cos\theta+\beta\sin\theta)-h_-$.
But $\alpha\cos\theta+\beta\sin\theta$ is the support function of
the point $(\alpha,\beta)$.  So
$(h_++\alpha\cos\theta+\beta\sin\theta)$ is a support function,
and $f$ is a difference of support functions as required.  Letting
$h_1=h_++\alpha\cos\theta+\beta\sin\theta$ and $h_2=h_-$ then
$f=h_1-h_2$ and for $i=1,2$ and by \eq{mu-bds}
$\|h_i''+h_i\|_{TV}=\|\mu_\pm\|\le 3\|f''+f\|_{TV}$.
\end{proof}

%%%%%%%%%%%%%%%%%%%%%%%%%%%%%%%%%%%%%%%%%%%%%%%%%%%%%%%%%%%%%%%%%%%%%%
\section{The affine isoperimetric inequality for arbitrary planar
convex sets.}\label{sec:aff-iso}
%%%%%%%%%%%%%%%%%%%%%%%%%%%%%%%%%%%%%%%%%%%%%%%%%%%%%%%%%%%%%%%%%%%%%%

If $h$ is a support function, then $h$ is Lipschitz and therefore
absolutely continuous.  Therefore the distributional derivative $h'$
of $h$ is just the classical derivative which exists almost
everywhere.  As $h$ is a support function then by
Proposition~\ref{spt-char} the second distributional derivative $h''$
is a measure and therefore $h'$ is of bounded variation.  By a
theorem of Lebesgue, the function $h'$ will be differentiable (in the
classical sense) almost everywhere.  Denote this derivative of $h'$ by
$Dh'$ to distinguish it from the distributional derivative.  In what
follows we will denote classical derivatives of a function $f$ by
$Df$.  As the first distributional derivative of $h$ agrees with the
classical derivative we have $Dh'=D^2h$ so that $Dh'$ is the second
classical derivative.

Recall, by a theorem of Alexandrov, a convex function on an $n$
dimensional space, and thus a support function, has a generalized
second derivative, called the Alexandrov second derivative, almost
everywhere and in the one dimensional case the Alexandrov second
derivative is just $D^2h$.

Various
authors~\cite{Leichtweiss:aff-area,Schutt-Werner:affine,Lutwak:aff-area,Werner:ill-affine}
have extended the definition of affine arclength (and more generally
higher dimensional affine surface area) from convex sets with $C^2$
boundary to general convex sets.  It was eventually shown all these
definitions are equivalent see i.e.~\cite{Dolzmann-Hug} and, for two
dimensional convex sets, are given in terms of the support function
%are given
by
$$
\int_\s (D^2h+h)^{2/3}\,d\theta.
$$
The following is the general form of the affine isoperimetric
inequality in the plane.

\begin{thm}\label{aff-iso}
Let $K$ be an compact convex body in the plane with area $A$ and
affine perimeter $\Omega$.  Then
\begin{equation}\label{eq:aff-iso}
\Omega\le 8\pi^2 A
\end{equation}
with equality if and only if $K$ is an ellipse.
\end{thm}

The proof of this Theorem is based on our Theorem~\ref{thm:main}
and the following result which compares the distributional and
classical derivatives of a support function.

\begin{pro}\label{decomp}
Let $h\cn\s\to \R$ be the support function of a bounded convex set.
Then distribution $h''+h$ is of the form
\begin{equation}\label{D2h}
h''+h=(D^2h+h)\,d\theta +d\mu
\end{equation}
where $d\theta$ is Lebesgue measure, the function $D^2h+h$ is in
$L^1(\s)$ and $\mu$ is a non-negative measure that is singular with
respect to Lebesgue measure (i.e. there is a set $N$ of Lebesgue
measure zero with $\mu(S\setminus N)=0$).
\end{pro}

The proof of the Proposition is elementary and is left to the
reader.

\begin{proof}[Proof of Theorem~\ref{aff-iso}]
Let $h_o$ be the support function of a planar convex body $K$. In
the proof of Lemma~\ref{lem:gen-F-h} we will see that it is
possible to choose $a_o$ and $b_o$ so that $h(\theta):=
h_o(\theta)+a_o\cos\theta+b_o\sin\theta$ is positive on $\s$ and
\begin{equation}\label{h2ortho}
\int_\s \frac{\cos\theta}{h^3(\theta)}\,d\theta=\int_\s
\frac{\sin\theta}{h^3(\theta)}\,d\theta =0.
\end{equation}
%(The idea is to choose $(a_o,b_o)$ to minimize
%$f(a,b)=\int_\s(h_o(\theta)+a\cos\theta+b\sin\theta)^{-2}\,d\theta$
%over the set of $(a,b)\in \R^2$ such that the origin is in the
%interior of $K+(a,b)$.)
Using the relation between $D^2h+h$ and $h''+h$ given by
Proposition~\ref{decomp} we have
$$
\int_\s h(D^2h+h)\,d\theta\le \int_\s h(h''+h)\,d\theta.
$$
This observation, preceded by H\"older's inequality,  gives,
%\begin{align*}
%\int_\s(D^2h+h)^{2/3}\,d\theta&=\int_\s h^{-2/3}h^{2/3}(D^2h+h)^{2/3}\,d\theta\\
%&\le \left(\int_\s\frac{d\theta}{h^2}\right)^{1/3}
    %\left(\int_\s h(D^2h+h)\,d\theta\right)^{2/3}\\
%&\le \left(\int_\s\frac{d\theta}{h^2}\right)^{1/3}
    %\left(\int_\s h(h''+h)\,d\theta\right)^{2/3}.
%\end{align*}
%Thus
\begin{equation}\label{ai1}
\left(\int_\s(D^2h+h)^{2/3}\,d\theta\right)^3 \le
\left(\int_\s\frac{d\theta}{h^2}\right)
    \left(\int_\s h(h''+h)\,d\theta\right)^{2}.
\end{equation}

In Theorem~\ref{thm:main} take $F=h^{-3}$.  Then \eq{h2ortho}
shows that conditions~\eqref{1.0} are satisfied.  Therefore
\begin{align}\nonumber
\left(\int_\s\frac{d\theta}{h^2}\right)^3
    \left(\int_\s h(h''+h)\,d\theta\right)
&=\left(\int_\s\frac{d\theta}{h^2}\right)^3\left(\int_\s
[h^2-(h')^2]\,d\theta\right)\\
\label{ai2}
&\le 4\pi^2\left(\int_\s \frac{d\theta}{h^2}\right)^2.
\end{align}
Combining \eq{ai1} and \eq{ai2} and using the fact that
$D^2h+h=D^2h_o+h_o$ and $\int_\s h(h''+h)\,d\theta= \int_\s
h_o(h_o''+h_o)\,d\theta$ gives
$$
\left(\int_\s(D^2h_o+h_o)^{2/3}\,d\theta\right)^3\le 4\pi^2
\int_\s h_o(h_o''+h_o)\,d\theta.
$$
This is the affine isoperimetric inequality for $K$.

If equality holds, then the equality conditions of
Theorem~\ref{thm:main} imply $h$ is the support function of an
ellipse centered at the origin.  Thus
$h_o=h-a\cos\theta-b\sin\theta$ is the support function of an
ellipse centered at $(-a,-b)$.
\end{proof}

%%%%%%%%%%%%%%%%%%%%%%%%%%%%%%%%%%%%%%%%%%%%%%%%%%%%%%%%%%%%%%%%%%%%%%
\section{A family of transforms.}\label{sec:trans}
%%%%%%%%%%%%%%%%%%%%%%%%%%%%%%%%%%%%%%%%%%%%%%%%%%%%%%%%%%%%%%%%%%%%%%

Let $\s$ be the unit circle in $\R^2$ with coordinate $\theta$ as
above.  For each $\lambda \in (0, \infty)$, let
$$
\psi_{\lambda}(\theta)= \sqrt{ \lambda^2 \cos^2 \theta +
\frac{1}{\lambda^2}\sin^2\theta} .
$$
Define on $\s$ a family of
mappings
$$
m_{\lambda}(\theta) =
\int_0^{\theta} \frac{dt}{\psi_{\lambda}^2 (t)} .
$$
When $\lambda = 1$, this is the identity map.
For $0 \leq \theta < \frac{\pi}{2}$
it is easy to verify that
\begin{equation}
m_{\lambda}(\theta) = \arctan \left( \frac{1}{\lambda^2} \tan \theta \right).
\label{2.7}
\end{equation}
For any measurable function $u$ on $\s$, define the transform
$$
(T_{\lambda} u)(\theta) = u(m_{\lambda}(\theta))
\psi_{\lambda}(\theta) .
$$

\begin{lem}\label{lem:m-T}  Let $u$ and $v$ be measurable functions on
$\s$ for which the integrals below exist.
%\mnote{\ralph Added a change suggested by Erwin?}
Then

(i) The mappings $m_{\lambda}(\cdot)$ each leave four points fixed:
$$
m_{\lambda}(0) = 0,\quad  m_{\lambda}\left(\frac{\pi}{2}\right) =
\frac{\pi}{2}, \quad m_{\lambda}(\pi) = \pi, \quad
m_{\lambda}\left(\frac{3\pi}{2}\right) = \frac{3\pi}{2} .
$$

(ii) The transforms leave the following integrals
invariant:
\begin{equation}
\int \frac{d\theta}{(T_{\lambda} u)^2} = \int \frac{d\theta}{u^2}\; ,
\label{2.1}
\end{equation}
\begin{equation}\label{2.a}
\int \frac{T_{\lambda} u}{(T_{\lambda}v)^3}\,d\theta = \int
\frac{u}{v^3}\,d\theta,
\end{equation}
\begin{equation}
\int_\s \{ (T_{\lambda}u)^2 - [(T_{\lambda}u)']^2 \}\,d\theta =
\int_\s [ u^2 - (u')^2 ] \,d\theta , \label{2.2}
\end{equation}
\begin{equation}
\int \frac{\cos \theta}{(T_{\lambda}u)^3(\theta)}  \,d\theta=
\frac{1}{\lambda} \int \frac{\cos \theta}{u^3(\theta)}   \,d\theta,
\label{2.3}
\end{equation}
\begin{equation}
\int \frac{\sin \theta}{(T_{\lambda}u)^3(\theta)}  \,d\theta=
\lambda \int \frac{\sin \theta}{u^3(\theta)}  \,d\theta .
\label{2.4}
\end{equation}
Here  ``$\int$'' represents the integral
with respect to $d\theta$ on any of the intervals $[0,
\frac{\pi}{2}]$, $[\frac{\pi}{2}, \pi],$ $ [\pi, \frac{3\pi}{2}],$ $
[\frac{3\pi}{2}, 2\pi],$ or $[0, 2\pi]$.  For \eq{2.1}, \eq{2.a},
\eq{2.3}, and \eq{2.4} $u$ and $v$ can be any measurable functions for
which the integrals converge.
In \eq{2.2} $u\in H^1(\s)$.

%In particular, if we set $\phi_{\lambda}(\theta) =
%\sqrt{\frac{1}{\lambda^2} \cos^2 \theta + \lambda^2 \sin^2 \theta
%}$, then
%$$ T_{\lambda}(\phi_{\lambda}(\theta) \equiv 1 . $$
%\label{lem2.1}
\end{lem}

\begin{proof}

(i) From \eq{2.7}, one can see that $m_\lambda(0)=0$,
$m_{\lambda}(\frac{\pi}{2}) = \frac{\pi}{2}$. Since the integrand is
symmetric about $\theta = \frac{\pi}{2}$ on $[0, \pi]$, and is
$\pi$-periodic, in follows that $\pi$ and $\frac{3\pi}{2}$ are also
fixed points of $m_\lambda$.

(ii) We verify invariance for the integrals on the interval $[0,
\frac{\pi}{2}]$. Then by the results in (i) and the symmetry of
$\psi_{\lambda}(\theta), \cos \theta,$ and $\sin \theta$, the
invariance of the integrals on the other intervals follows.
\smallskip

Equations \eq{2.1} and \eq{2.a} are direct consequences from the
substitution $\tilde{\theta} = m_{\lambda}(\theta)$.

As $C^2(\s)$ is dense in $H^1(\s)$ it is enough to verify \eq{2.2} in
the case $u\in C^2(\s)$.  After integrating by parts, we only need to
show
$$
\int_\s T_{\lambda} u [ T_{\lambda} u + (T_{\lambda} u)'' ]\,d\theta
= \int_\s u ( u + u'' )\,d\theta .
$$
We employ the fact that $\psi_{\lambda}$ is a solution of the
equation
\begin{equation}\label{2.9}
\psi_{\lambda}'' + \psi_{\lambda} = \frac{1}{\psi_{\lambda}^3} .
\end{equation}
It follows from \eq{2.9} and a straightforward calculation
that
$$
T_{\lambda} u [ T_{\lambda} u + (T_{\lambda} u)'' ] =
u(m_{\lambda}(\theta)) [ u(m_{\lambda}(\theta)) +
u''(m_{\lambda}(\theta)) ] \frac{1}{\psi_{\lambda}^2(\theta)} .
$$
Again using the change of variable $\tilde{\theta} =
m_{\lambda}(\theta)$, we see \eq{2.2} holds.

To obtain \eq{2.3} and \eq{2.4}, we write $ \tan \theta = \lambda^2
\tan \tilde{\theta}$. It follows that
$$
\frac{\cos \theta}{\psi_{\lambda}(\theta)}
 = \frac{1}{\sqrt{\lambda^2 + \frac{1}{\lambda^2} \tan^2 \theta}} =
\frac{1}{\lambda} \cos \tilde{\theta},\quad
\frac{\sin\theta}{\psi_\lambda(\theta)}
=\frac{1}{\sqrt{\frac{\lambda^2}{\tan^2\theta}+\frac1\lambda}}
=\lambda\sin\tilde\theta
$$
and another application of the substitution $\tilde{\theta} =
m_{\lambda}(\theta)$ completes the proofs of \eq{2.3} and
\eq{2.4}.
\end{proof}

\begin{rem}\label{rem:mass}
Let $u$ be a positive continuous function on $\s$, then the integral
of $\frac{d\theta} {(T_\lambda u)^2}$ is independent of $\lambda$.  It
is not hard to check that as $\lambda\to\infty$ the mass of
$\frac{d\theta} {(T_\lambda u)^2}$ concentrates about the points
$\pi/2$ and $3\pi/2$ and when $\lambda\to 0$ the mass concentrates
about $0$ and $\pi$.
\end{rem}

%%%%%%%%%%%%%%%%%%%%%%%%%%%%%%%%%%%%%%%%%%%%%%%%%%%%%%%%%%%%%%%%%%%%%%
\section{Proof of the main inequality.}
\label{sec:proof}
%%%%%%%%%%%%%%%%%%%%%%%%%%%%%%%%%%%%%%%%%%%%%%%%%%%%%%%%%%%%%%%%%%%%%%

\subsection{Some lemmas.}

\begin{lem}\label{lem2.1}
If $\{u_k\}$ is a bounded sequence in $H^1(\s)$, then there exists a
subsequence (still denoted by $\{u_k\}$) and $u_o\in H^1(\s)$, such
that $u_k\to u_o$ in the weak topology of $H^1(\s)$,
\begin{equation} \label{2.5}
u_k \to u_o \quad \text{ in } C^{\beta}(\s) \quad \text{for all}\quad \beta <
\frac{1}{2}.
\end{equation}
This implies
\begin{equation}\label{up-semi}
\limsup_{k\to\infty}\int_\s [u_k^2-(u_k')^2]\,d\theta\le
\int_\s [u_o^2-(u_o')^2]\,d\theta.
\end{equation}
Moreover, if $u_o(\theta_0) = 0$ at some point $\theta_0$, then for $\delta>0$
\begin{equation}\label{2.6}
\int_{\theta_0}^{\theta_0+\delta}\frac{d\theta}{u_o^2}=\infty
\quad\text{and}\quad
\int_\s \frac{d\theta}{u_k^2} \to \infty \quad \text{as}\quad k \to \infty.
\end{equation}
\end{lem}

\begin{lem}\label{3z}
Assume $\{u_k\}$ is a bounded sequence in $H^1(\s)$ with $u_k>0$,
$u_k\to u_o\in H^1(\s)$ in the weak topology,
\begin{equation}\label{uk-mass}
\int_\s \frac{\cos\theta}{u_k^3}\,d\theta=0
=\int_\s \frac{\sin\theta}{u_k^3}\,d\theta,
\end{equation}
and that $u_o$ has at least one zero.  Then, viewing the zeros of
$u_o$ as a subset of $\s\subset \R^2$,
\begin{equation}\label{co-h}
(0,0)\in \text{convex hull of the zeros of $u_o$.}
\end{equation}
If $u_o$ has three or more zeros then
\begin{equation}\label{neg}
\int_{\s} [ u_k^2 - (u_k')^2 ]\, d\theta  < 0 \quad \text{for
sufficiently large $k$}.
\end{equation}
\end{lem}

\begin{lem}\label{lem:gen-F-h}
Suppose the inequality \eq{1.2} holds under the
stronger conditions:

i) $F$ is measurable and positive on $\s$ and $h \in H^1(\s)$ is
positive;
%\mnote{\ralph A bit or rewriting here.}

ii) $F$ satisfies the orthogonality conditions \eq{1.0}
and $h$ satisfies orthogonality conditions
\begin{equation}
\int_\s \frac{\cos \theta}{h^3}  \, d\theta = 0 = \int_\s
\frac{\sin \theta}{h^3}  \, d\theta .
\label{h}
\end{equation}

Then the same inequality \eq{1.2} holds without the
orthogonality conditions \eq{h} on $h$ and the strict positivity of
$F$.
%\mnote{\ralph Made a change suggested by Gaoyong. Old version still in file.}
%, and  only assuming that
%$F$ is non-negative and $h\in H^1(\s)$ is non-negative.
\end{lem}

\begin{proof}[Proof of Lemma \ref{lem2.1}.]
That there is a $u_o\in H^1(\s)$ and a subsequence with $u_k\to u_o$
in the weak topology follows from the weak compactness of the closed
balls in a Hilbert space.  Then \eq{2.5} is a direct consequence of
the compact Sobolev imbedding of $H^1(\s)$ into $C^{\beta}(\s)$ for
any $\beta < \frac{1}{2}$. To prove\eq{up-semi} use the fact that the
norm of a Hilbert space is lower semi-continuous with respect to weak
convergence and thus $\liminf_{k\to\infty}\int_\s [u_k']^2\,d\theta\ge
\int_\s [u_o']^2\,d\theta$.  From\eq{2.5} $\lim_{k\to\infty}\int_\s
u_k^2\,d\theta=\int_\s u_o^2\,d\theta$.  Together these
imply\eq{up-semi}.

Assume that $u_o$ vanishes at $\theta_0$.  Then by the Sobolev
imbedding $H^1(\s)\subset C^{\frac12}(\s)$, or an elementary H\"older
inequality argument,  $u_o\in C^{\frac12}(\s)$ and therefore $ |
u_o(\theta) |=|u_o(\theta)-u_o(\theta_0)| \leq C_1
\sqrt{|\theta-\theta_0|} $ which implies the divergence of the
integral $\int_{\theta_0}^{\theta_0+\delta}u_o^{-2} {d\theta}$.
As $u_k\to u_o$ uniformly this implies the second part of \eq{2.6} and
completes the proof of the Lemma.
\end{proof}

\begin{proof}[Proof of Lemma \ref{3z}.]
Because the imbedding of $H^1(\s)$ into $C^\beta(\s)$ is compact for
$\beta\in[0,1/2)$ the weak convergence $u_k\to u_o$ implies $\{u_k\}$
converges to $u_o$ uniformly.  By Lemma~\ref{lem2.1} the integral
$\int_\s u_o^{-2}\,d\theta$ diverges and therefore
$\int_\s u_o^{-3}\,d\theta$ also diverges.  Thus
$\int_\s u_k^{-3}\,d\theta\to \infty$ as $k\to\infty$.  Let
$c_k:=\left(\int_\s u_k^{-3}\,d\theta\right)^{-1}$.  Then
$c_k{u_k^{-3}(\theta)}{d\theta}$ is a probability measure on $\s$
and the conditions \eq{uk-mass} imply the center of mass of this
measure is $(0,0)$.  But as $k\to\infty$ the masses of the measures
$c_k{u_k^{-3}(\theta)}{d\theta}$ concentrate at the zeros of
$u_0$. This implies \eq{co-h}.

If $u_o$ has three or more zeros, then the convex hull property
\eq{co-h} implies there are three zeros $\theta_1$,
$\theta_2$, $\theta_3$ of $u_o$ such that
\begin{equation}
\text{The zeros $\theta_1,\theta_2,\theta_3$ of $u_o$ are not on an
arc of length less than $\pi$.}
\label{3.B}
\end{equation}
We will show this implies
\begin{equation}\label{3.7}
\int_{\s} [ u_o^2 - (u_o')^2 ] \, d\theta < 0 .
\end{equation}
Which, by \eq{up-semi} of Lemma~\ref{lem2.1}, implies \eq{neg}.

To see \eq{3.7}, we write the integral in three parts:
\begin{equation}
\int_{\s} [ u_o^2 - (u_o')^2 ]\, d\theta  = \left\{ \int_{\theta_1}^{\theta_2} +
\int_{\theta_2}^{\theta_3} + \int_{\theta_3}^{\theta_1} \right\} [ u_o^2
- (u_o')^2 ]\, d\theta  = I_1 + I_2 + I_3 . \label{3.8}
\end{equation}
From \eq{3.B}, we see that lengths of intervals of integration in
\eq{3.8} are all less than or equal to $\pi$, and at least two of them
are strictly less than $\pi$.
%\mnote{\ralph Changed the eigenvalue
%argument to the reference to \cite{Hardy-Littlewood-Polya}.  Old version still in file.}
But, \cite[p.~185]{Hardy-Littlewood-Polya},
if $\theta_{i+1}-\theta_i\le\pi$, then
$u_o(\theta_i)=u_o(\theta_{i+1})=0$ implies
$\int_{\theta_i}^{\theta_{i+1}} [ u_o^2 - (u_o')^2 ] \, d\theta<0 $
unless $\theta_{i+1}-\theta_i=\pi$ and $u_o=C\sin\theta$ on
$[\theta_i,\theta_{i+1}]$. This proves \eq{3.7} and completes the
proof of the lemma.
%Combining this with the fact
%that the smallest eigenvalue of the eigenvalue problem
%\begin{equation}
%\left\{\begin{array}{ll} - u'' = \lambda u  & \text{on}\  [a, b] \\
%u(a) = 0 = u(b) \end{array} \right. \label{3.9}
%\end{equation}
%is $\lambda_1 = \left(\pi/(b-a)\right)^2$; we conclude that
%$I_1$, $I_2$, and $I_3$ are all non-positive, and at least two of them
%are negative.  This proves \eq{3.7} and completes the proof of the
%lemma.
\end{proof}

\begin{proof}[Proof of Lemma~\ref{lem:gen-F-h}.]
We assume that the inequality \eq{1.2} holds under the assumptions i)
and ii) of Lemma~\ref{lem:gen-F-h}.  We first claim that for each
positive function $h \in C^2(\s)$, there exists an $h_o(\theta) = a_o
\cos\theta + b_o \sin\theta + h(\theta)$ that satisfies the orthogonality
conditions \eq{h}. To see this minimize the function
$$
f(a, b) = \int_\s \frac{1}{( a \cos\theta + b \sin\theta + h(\theta)
)^2} d\theta
$$
for any real numbers $a$ and $b$, such that $a \cos\theta + b
\sin\theta + h(\theta) > 0$ for all $\theta$.  It is obvious that
$f(a, b)$ is bounded from below by zero. Let $\{h_k(\theta) = a_k
\cos\theta + b_k \sin\theta + h(\theta)\}$ be a minimizing
sequence. From $h_k(t) > 0$, one can easily see that $\{a_k\}$ and
$\{b_k\}$ are bounded, and hence there exist subsequences converging
to some $a_o,b_o\in\R$.  Then $(a_o, b_o)$ is a minimizer of $f$.

Moreover, from Lemma \ref{lem2.1}, we can see that $h_o(\theta) = a_o
\cos\theta + b_o \sin\theta + h(\theta) > 0$. (Otherwise $h_o$ has a
zero and (by Lemma~\ref{lem2.1}) $\int_\s h_o^{-2}{d\theta}=\infty$,
contradicting that $h_o$ is a minimizer.)  Consequently, at $(a_o,
b_o)$, we have ${\partial f}/{\partial a} = 0 = {\partial
f}/{\partial b}$. This implies the orthogonality conditions \eq{h} on
$h_o$.

We now show that if inequality \eq{1.2} holds for $h_o =
a_o \cos\theta + b_o \sin\theta + h(\theta)$, then it is also holds for $h$.
By the orthogonality conditions \eq{1.0}
on $F$.
\begin{equation}
\int_\s F(\theta)h_o(\theta)\, d\theta = \int_\s
F(\theta)h(\theta)\, d\theta ,
\label{5.2}
\end{equation}
and if $h$ is of class $C^2$ we can use  use integration by parts and the fact
that both $\sin \theta$ and $\cos \theta$ are in the kernel of the differential
operator ${d^2}/{d\theta^2} + 1$ to get
\begin{align}
\int_\s [h_o^2-(h_o')^2]\,d\theta
  &=\int_\s h_o( h_o'' + h_o )\, d\theta\nonumber\\
& = \int_\s h( h'' + h )\, d\theta =\int_\s [h^2-(h')^2]\,d\theta.
\label{5.3}
\end{align}
This will also hold for $h\in H^1(\s)$ by approximating by $C^2$
functions.  So if \eq{1.2} holds for $h_o$ and $F$, then \eq{5.2} and
\eq{5.3} show it holds for $h$ and $F$.

To see that inequality \eq{1.2} holds also for non-negative
continuous functions $F$ and non-negative $h\in H^1(\s)$, we let
$$
F_{\epsilon} = F + \epsilon \;\;\; \text{ and } \;\; h_{\epsilon}
= h + \epsilon .
$$
Then obviously, for each $\epsilon > 0$, both $F_{\epsilon}$ and
$h_{\epsilon}$ are positive, and $F_{\epsilon}$ satisfies the
orthogonality conditions \eq{1.0}. Therefore inequality
\eq{1.2} holds for $F_{\epsilon}$ and $h_{\epsilon}$. Take the
limit as $\epsilon\to0$ to see that see that \eq{1.2} is also holds
for $F$ and $h$.

Finally the extensions to $F$ non-negative and measurable follows by
approximating $F$ by positive functions satisfying the orthogonality
conditions\eq{1.0} and taking limits.
\end{proof}

\subsection{Outline of the Proof.}

Let
$$
G = \left\{ v : v > 0, \; \int_S \frac{1}{v^3} d\theta < \infty ,  \; \int_{S} \frac{1}{v^3}
\cos\theta \, d\theta = 0 = \int_{S} \frac{1}{v^3} \sin\theta \,
d\theta \right\}$$ 
and let
$$
\tilde{G} = \{ u \in G : u \in H^1(S) \}.
$$

Define on $\tilde{G} \times G$,
$$I(u, v) = \dfrac{\{ \int_{S}  \frac{1}{v^2} d\theta \}^3 \{
\int_{S} [ u^2 - (u')^2 ] d\theta \}}{\{ \int_{S} \frac{u}{v^3}
d\theta \}^2} .$$ It is obvious that for any constants $t$ and
$s$, we have
$$ I(tu, sv) = I(u, v).$$
To prove the theorem, it is equivalent to show that
\begin{equation}
I(u, v) \leq 4 \pi^2 \;\;\; \forall (u, v) \in \tilde{G} \times G.
\label{1.1}
\end{equation}
and the equality holds if and only if
\begin{equation}
u = k_1 \sqrt{\lambda^2 \cos^2\theta +
\frac{1}{\lambda^2}\sin^2\theta }
\label{1.a}
\end{equation}
and
\begin{equation}
v = k_2 \sqrt{\lambda^2 \cos^2\theta +
\frac{1}{\lambda^2}\sin^2\theta }
\label{1.b}
\end{equation}
with any non-negative constants $k_1$, $k_2$,  and $\lambda > 0$.
\bigskip

First,  we show that there exists a constant $C < \infty$, such
that
$$ I(u, v) \leq C \;\;\;  \forall (u, v) \in \tilde{G} \times G .$$
It is done by applying the family of transforms and by using a
contradiction argument.

Then we study a maximizing sequence $\{(u_k, v_k)\}$ of the
functional $I(u, v)$. Usually, such a sequence may be unbounded.
However thanks to the family of transforms, we are able to convert
it into a new sequence which converges to a maximum $(u_o, v_o)$
in $\tilde{G} \times G$. Finally, we use the well-known
classification results on the solutions of the corresponding
Euler-Lagrange equations to arrive at the conclusion of the
theorem.
\bigskip

\subsection{The Proof}

{\bf Part I.}

In this part, we show that there exists a constant $C < \infty$,
such that
\begin{equation}
I(u, v) \leq C , \;\; \forall (u, v)  \in \tilde{G} \times G.
\label{3.1}
\end{equation}

We argue by contradiction. Suppose in contrary, there exists a
sequence $\{(\tilde{u}_k, v_k)\}$ in $\tilde{G} \times G$, such
that $I(\tilde{u}_k, v_k) \ra \infty$, as $k \ra \infty$.

Let $u_k = \dfrac{\tilde{u}_k}{\|\tilde{u}_k\|}$, then
\begin{equation}
\|u_k\| = 1 , \;\; \mbox{ and } I(u_k, v_k) \ra \infty .
\label{3.2}
\end{equation}

Here and in the rest of the paper, for convenience of writing, we
use $\{(u_k, v_k)\}$ to denote the sequence itself or one of its
subsequences.

 The first part of (\ref{3.2}) implies that
 \begin{equation}
 J(u_k) = \int_\s [ u_k^2 -(u_k')^2 ]\,d\theta \mbox{ is bounded },
 \label{3.A}
 \end{equation}
 Therefore from the second part, we must
have
\begin{equation}
u_k(\theta) \ra 0 ,\;\; \mbox{ for some } \theta.
\label{3.4}
\end{equation}

Otherwise, if $\{u_k\}$ is bounded away  from zero, then by the
H\"older inequality
\begin{equation}
 \int_\s \frac{1}{v^2}\, d\theta  \leq \left( \int_\s \frac{u}{v^3}\, d\theta \right)^{2/3}
\left(\int_\s \frac{1}{u^2}\, d\theta \right)^{1/3} \label{A} 
\end{equation} 
we would arrive at the boundedness of
$$
\left( \int_\s \frac{1}{v_k^2}\,d\theta \right)^3 
\left( \int_\s \frac{u_k}{v_k^3}\,d\theta  \right)^{-2} .
$$
This, together with (\ref{3.A}), contradicts with the second part
of (\ref{3.2}).

Without loss of generality we assume that $J(u_k)>0$ for all $k$.
We have shown that $u_o$ has at least one zero and therefore by
the convex hull property of Lemma~\ref{3z} the point $(0,0)$ is in
the convex hull of the zeros of $u_o$.  This implies that $u_0$
has at least two zeros.  If $u_o$ has three or more zeros then
\eq{neg} of Lemma~\ref{3z} implies that $J(u_k)<0$ which is not
the case.  Thus $u_o$ has exactly two zeros.

As $u_o$ has exactly two zeros, the convex hull property \eq{co-h}
implies the two zeros must be antipodal, say they are at $\theta =
\frac{\pi}{2}$ and $\frac{3\pi}{2}$. Obviously, at these two
points, $u_k^{-2} \to \infty$. For each $u_k$, pick a point $p_k$
near ${\pi}/{2}$, such that
\begin{equation}
\int_{p_k-\frac{\pi}{2}}^{p_k} \frac{1}{u_k^2}\, d\theta=
\int_{p_k}^{p_k+\frac{\pi}{2}} \frac{1}{u_k^2}\,d\theta.
\label{3.10}
\end{equation}
Then, $p_k \to {\pi}/{2}.$ For a number $\delta > 0$ (to be
concrete $\delta=\pi/4$ will work), let
$$
D_k^1 = \left\{ \theta : p_k - \frac{\pi}{2} \leq \theta \leq
p_k - \delta \right\} , \quad D_k^2 = \left\{ \theta : p_k +
\delta \leq \theta \leq p_k + \frac{\pi}{2} \right\}
$$
$$
D_k = D_k^1 \cup D_k^2 .
$$
We will apply the family of transforms $T_\lambda$ introduced in
Section~\ref{sec:trans}. We say that the family of transforms
$T_{\lambda}$ in Lemma~\ref{lem2.1} are centered at
$\frac{\pi}{2}$, and write $T_{\lambda} = T_{\lambda,
\frac{\pi}{2}}$. Similarly, one can define transforms centered at
any point $q$, and denote them by $T_{\lambda, q}$.

By Lemma~\ref{lem:m-T} and Remark~\ref{rem:mass} for each $u_k$,
one can choose a transform $T_k = T_{\lambda_k, p_k}$, such that
\begin{equation}\label{3.11}
\int_{B_{\delta}(p_k)} \frac{d\theta}{(T_k u_k)^2} = \int_{D_k}
\frac{d\theta}{(T_k u_k)^2} ,
\end{equation}
where $B_{\delta}(p_k) = ( p_k - \delta, p_k + \delta ).$ Let $w_k
= \|T_k u_k\|_{H^1}^{-1}\,T_k u_k$. Then $\|w_k\|_{H^1} = 1$ and
we can apply Lemma~\ref{lem2.1} to the sequence $\{w_k\}$ and find
a subsequence, still denoted by $\{w_k\}$, and a $w_o\in H^1(\s)$
such that $ w_k \to w_o $ in the weak topology of $H^1$ and $
w_k(\theta) \to w_o(\theta)$ in $ C^{\beta}$ for all $\beta <
\frac{1}{2}$.
\medskip

If $w_o$ has no zeros, then by (\ref{2.1}) in Lemma \ref{lem:m-T}
and (\ref{A}), $I(u_k, v_k)$ is bounded, and we are done.

Therefore, we may assume that $w_o$ has at least one zero. By the
convex hull property of Lemma~\ref{3z}
\begin{equation}\label{co-h-w}
(0,0)\in \text{convex hull of the zeros of $w_o$.}
\end{equation}
This implies that $w_o$ has at least two zeros and if $w_o$ has
three or more zeros then Lemma~\ref{3z} implies $J(w_k)<0$ for
large $k$, which, by (\ref{2.2}) of Lemma \ref{lem:m-T},  again
would contradict with the assumption that $J(u_k) > 0$. Therefore
$w_o$ has exactly two zeros and by the convex hull property
\eq{co-h-w} these zeros are antipodal. Let the zeros be $\theta_0$
and $\theta_1$ and we can assume that $\theta_0\in [0,\pi]$. Then
$\int_0^\pi w_o^{-2}\,d\theta=\infty$ (by \eq{2.6}), $w_k\to w_o$
uniformly, and $p_k\to\pi/2$ imply
\begin{equation}\label{3.13}
\int_{p_k -\frac{\pi}{2}}^{p_k +\frac{\pi}{2}} \frac{
d\theta}{w_k^2} \to \infty .
\end{equation}
From \eq{3.10}, \eq{3.11} and the properties of the transforms, we
also have
\begin{equation}
\int_{p_k -\frac{\pi}{2}}^{p_k} \frac{d\theta}{w_k^2}  =
\int_{p_k}^{p_k +\frac{\pi}{2}} \frac{d\theta}{w_k^2} ,
\label{3.14}
\end{equation}
and
\begin{equation}
\int_{B_{\delta}(p_k)} \frac{d\theta}{w_k^2}  = \int_{D_k}
\frac{d\theta}{w_k^2}  . \label{3.15}
\end{equation}
Now \eq{3.13}, \eq{3.14}, and \eq{3.15} imply that the integrals
of $w_k^{-2}$ on all the four sets
$$
[p_k-\frac{\pi}{2}, p_k],\  [p_k, p_k +\frac{\pi}{2}],\
B_{\delta}(p_k), \text{ and } D_k
$$
approach infinity. Therefore $w_o$ has at least one zero on each
of the following sets
\begin{equation}\label{zero-bds}
[0, \frac{\pi}{2}],\  [\frac{\pi}{2}, \pi],\
B_{\delta}(\frac{\pi}{2}), \text{ and } [0, \pi] \setminus
B_{\delta}(\frac{\pi}{2}).
\end{equation}
From this we see that $w_o$ has at least two zeros on the closed
upper half circle.  As the two zeros of $w_o$  are antipodal they
must be $0$ and $\pi$.  But as $B_{\delta}({\pi}/{2})$ also
contains a zero this implies that $w_o$ has three zeros, a
contradiction.
\bigskip

{\bf Part II.}

In part I, we have shown that
$$I(u, v) \leq C < \infty \;\; \forall \; (u, v) \in \tilde{G} \times  G. $$
To obtain the least possible value of the constant C, we consider
a maximizing sequence $\{(\tilde{u}_k, v_k)\}$ with
$\|\tilde{u_k}\|= 1.$

Using an entire similar argument as in part I, we can show that
there exists a family of transforms $T_k = T_{\lambda_k, p_k}$,
such that for $u_k = \dfrac{T_k \tilde{u_k}}{\|T_k
\tilde{u_k}\|}$, we have
$$
u_k(\theta) \ra u_o(\theta) > 0 ,  \;\; \forall
\theta \in [0, 2\pi]. 
$$

It follows from (\ref{A}) that
\begin{equation}
I(u_k, v_k) \leq I(u_k, u_k) .
\label{4.b}
\end{equation}
Thus, $(u_o, u_o)$ is a maximum of $I(u, v)$;  and
furthermore, it is in the interior of $\tilde{G} \times
\tilde{G}$. Therefore, we have
$$
\frac{\partial I}{\partial u} \bigg|_{\tilde{G}\times \tilde{G}} (u_o,
u_o) = 0.
$$
Through a straightforward calculation, one can see that a
constant multiple of $u_o$, still denoted by $u_o$, satisfies the
following Euler-Lagrange equation

\begin{equation}
u_o''(\theta) + u_o(\theta) = \frac{1}{u_o^3} + a \frac{\cos
\theta}{u_o^4} + b \frac{\sin \theta}{u_o^4} . \label{4.1}
\end{equation}

To determine the constants a and b, we multiply both sides of
(\ref{4.1}) by $\cos \theta$ and $\sin \theta$ respectively, then
integrate over $[0, 2\pi]$ to obtain
\begin{equation}
a \int_\s \frac{\cos^2 \theta}{u_o^4}\,d\theta + b \int_\s \frac{\sin \theta \cos
\theta}{u_o^4}\,d\theta = 0 , \label{4.2}
\end{equation}
\begin{equation}
a \int_\s \frac{\cos \theta \sin \theta}{u_o^4}\,d\theta + b \int_\s \frac{\sin^2
\theta}{u_o^4}\,d\theta = 0 . \label{4.3}
\end{equation}

Using the H\"older inequality, one can show that
\begin{equation}
\left|\begin{array}{cc}
\int_\s \frac{\cos^2 \theta}{u_o^4}\,d\theta &
\int_\s \frac{\sin \theta \cos \theta}{u_o^4}\,d\theta \\
\int_\s \frac{\sin \theta \cos \theta}{u_o^4}\,d\theta  & \int_\s \frac{\sin^2
\theta}{u_o^4}\,d\theta
\end{array} \right| > 0 .\label{4.4}
\end{equation}

Therefore the algebraic system (\ref{4.2}) and (\ref{4.3}) has
only the trivial solution $a = b = 0$. Consequently,  $u_o$ satisfies
\begin{equation}
u_o'' + u_o = \frac{1}{u_o^3} . \label{4.5}
\end{equation}

Now by the well-known classification result for equation
(\ref{4.5}), we have
$$ u_o(\theta) = \sqrt{\lambda^2 \cos^2 \theta +
\frac{1}{\lambda^2} \sin^2 \theta } $$ for some constant
$\lambda$.

Then, a straight forward calculation leads to
$$I(u_o, u_o) = 4 \pi^2 .$$

Finally, since in (\ref{A}), the equality holds if and only if v
is a constant multiple of u,  we see that if $(u_o, v_o)$ is a
maximum of the functional $I(u, v)$,  then $v_o$ must be a
constant multiple of $u_o$.

This completes the proof of the theorem.

\bigskip

{\bf Acknowledgments:} It is a pleasure to thank Steve Dilworth
and Anton Schep for various references. We also thank Anton and
Kostya Oskolkov for helpful comments regarding the exposition.

%\bibliographystyle{amsplain}
%\bibliography{/home/howard/tex/inputs/HowRefs.bib}

\providecommand{\bysame}{\leavevmode\hbox to3em{\hrulefill}\thinspace}
\providecommand{\MR}{\relax\ifhmode\unskip\space\fi MR }
\providecommand{\MRhref}[2]{%
  \href{http://www.ams.org/mathscinet-getitem?mr=#1}{#2}
}
\providecommand{\href}[2]{#2}

\end{document}